\documentclass[12pt]{article}
\usepackage{latexsym,a4wide,amssymb,amsmath,graphicx,epsfig}

\newtheorem{thm}{Theorem}
\newtheorem{lem}{Lemma}
\newtheorem{cor}{Corollary}
\newtheorem{conj}{Conjecture}
\newtheorem{prop}{Proposition}
\newtheorem{exer}{Exercise}
\newtheorem{defi}{Definition}
\newtheorem{exam}{Example}

\newcommand{\ebox}{\hfill $\Box$\\\vspace{0.15cm}}
\newcommand{\pr}{{\bf Proof.}\ }
\newcommand{\bt}{\begin{thm}}
\newcommand{\et}{\end{thm}}
\newcommand{\bl}{\begin{lem}}
\newcommand{\el}{\end{lem}}
\newcommand{\bp}{\begin{prop}}
\newcommand{\ep}{\end{prop}}
\newcommand{\bc}{\begin{cor}}
\newcommand{\ec}{\end{cor}}
\newcommand{\bcj}{\begin{conj}}
\newcommand{\ecj}{\end{conj}}
\newcommand{\bex}{\begin{exer}}
\newcommand{\eex}{\end{exer}}
\newcommand{\bi}{\begin{itemize}}
\newcommand{\ei}{\end{itemize}}
\newcommand{\be}{\begin{equation}}
\newcommand{\ee}{\end{equation}}
\newcommand{\ben}{\begin{enumerate}}
\newcommand{\een}{\end{enumerate}}
\newcommand{\bd}{\begin{defi}}
\newcommand{\ed}{\end{defi}}
\newcommand{\bexam}{\begin{exam}}
\newcommand{\eexam}{\end{exam}}

\newcommand{\mt}{t\kern-0.035cm\char39\kern-0.03cm}
\newcommand{\ml}{l\kern-0.035cm\char39\kern-0.03cm}
\newcommand{\md}{d\kern-0.035cm\char39\kern-0.03cm}

\newcommand{\noi}{\noindent}

\begin{document}

\title{\vspace{-2.3cm} A note on the superadditive and the subadditive transformations of aggregation functions}

\author{}
\date{}
\maketitle

\begin{center}
\vspace{-1.3cm}

{\large Alexandra \v{S}ipo\v{s}ov\'a \\
\vspace{2mm}

 {\small Slovak University of Technology}

{\small Faculty of Civil Engineering}

 {\small Department of Mathematics and Descriptive Geometry}

 {\small Radlinsk\'eho 11, 810 05 Bratislava}

{\small   Slovakia}

{\small alexandra.siposova@stuba.sk}
}
\vspace{4mm}

\end{center}

\begin{abstract}

We expand the theoretical background of the recently introduced superadditive and subadditive transformations of aggregation functions $A$. Necessary and sufficient conditions ensuring that a transformation of a proper aggregation function is again proper are deeply studied and exemplified. Relationships between these transformations are also studied.

\vskip 2mm

\noi {\em Keywords: aggregation function, subadditive transformation, superadditive transformation}

\end{abstract}

\vskip 3mm

\section{Introduction}\label{intro}

Motivated by applications in economics, subadditive and superadditive transformations of aggregation functions on $R^+=[0, \infty [$ have been recently introduced in ${\cite{Greco}}$.  Formally, both these transformations can be introduced on the improper real interval $[0, \infty]$.

\bd A mapping $A: [ 0, \infty]^n \to [0, \infty]$ is called an ($n$-ary) aggregation function if $A(0, \ldots , 0)=0$ and $A$ is increasing in each coordinate. Further, $A$ is called a proper ($n$-ary) aggregation function if it satisfies the following two additional constraints:
\medskip

${\rm (i)}$ \ \ $A({\bf x}) \in\ ]0, \infty [$ for some ${\bf x} \in\ ]0, \infty[^n$,
\smallskip

${\rm (ii)}$ \ \ $A(\bf x) < \infty $ for all ${\bf x} \in [0, \infty[^n$. \ed

Though for real applications we only need proper aggregation functions (in fact, their restriction to the domain $[0, \infty[^n$), a broader framework of all ($n$-ary) aggregation function is of advantage in a formal description of our results, making formulations and expressions more transparent. Observe that our framework is broader than the concept of aggregation functions on $[0, \infty]$ as introduced in ${\cite{Bel}}, {\cite{Grab}}$, which does not cover Sugeno integral based aggregation functions, for example. We denote the class of all $n$-ary aggregation function by ${\cal A}_n$, and the class of all $n$-ary proper aggregation function by ${\cal P}_n$.
\medskip

The next definition was motivated by optimization tasks treated in linear programming area and related areas ${\cite{Den}}$, as well as by recently introduced concepts of concave ${\cite{Leh}}$ and convex ${\cite{Mes}}$ integrals.

\bd For every $A \in {\cal A}_n$ the subadditive transformation $A_*:[0,\infty ]^n \to [0, \infty]$ of $A$ is given by

\begin{equation}
A_*({\bf x})= {\rm inf}\  \{\sum_{i=1}^k A({\bf y}^{(i)}) \ \  | \ \ \sum_{i=1}^k {\bf y}^{(i)} \geq {\bf x}\} \label{sub} \end{equation}

\noindent Similarly, for every $A\in {\cal A}_n$ the superadditive transformation $A^*:[0,\infty ]^n \to [0, \infty]$ of $A$ is defined by

\begin{equation}
A^*({\bf x})= {\rm sup}\  \{\sum_{j=1}^{\ell} A({\bf y}^{(j)}) \ \ | \ \ \sum_{j=1}^{\ell} {\bf y}^{(j)} \leq {\bf x}\}\ . \label{super}\end{equation}
\ed

Observe that the transformation (\ref{sub}) was originally introduced in \cite{Greco} for $A \in {\cal K}_*^n$, where ${\cal K}_*^n$ is the class of all $n$-ary proper aggregation functions (restricted to $[0, \infty[^n$) such that also $A_*$ is proper, that is, $A_* \in {\cal P}_n$. Similarly, $A^*$ given by  (\ref{super}) was originally introduced in \cite{Greco} only for $A \in {\cal K}_n^*$, where ${\cal K}_n^*$ is the class of all $A \in {\cal P}_n$ (restricted to $[0, \infty[^n$), so that $A^* \in {\cal P}_n$ as well.
\smallskip

Theorem $2$ in \cite{Greco} gives a necessary and sufficient condition ensuring that a function $A \in {\cal P}_n$ has also the property that $A \in {\cal K}_n^*$. We develop this result, giving an equivalent condition. Moreover, we also characterize all the functions $A \in {\cal P}_n$ such that $A \in {\cal K}_*^n$. Our approach is based on a deep study of transformations (\ref{sub}) and (\ref{super}) on unary aggregation functions that belong to ${\cal P}_1$. Our approach allows to show that for any $A \in {\cal P}_n$ we have the inequality $(A_*)^* \leq (A^*)_*$.
\smallskip

The paper is organized as follows. In the next section, the classes ${\cal K}_1^*$ and ${\cal K}_*^1$ are completely described, showing that the properties in a neighbourhood of $0$ are important for characterization of elements of these classes. In Section \ref{s3}, necessary and sufficient conditions for a function $A \in {\cal P}_n$ to belong to ${\cal K}_n^*$, or to ${\cal K}_*^n$, are given. Section \ref{s4} is devoted to the study of relationships of transformations $(A_*)^*$ and $(A^*)_*$. Finally, some concluding remarks are added.
\smallskip

\section{The one-dimensional case}\label{s2}

We begin with basic results which show how the values of the subadditive and superadditive transformations of one-dimensional aggregation functions depend on the behavior of the functions near zero.

\bt\label{ab} Let $h$ be an unary aggregation function on $[ 0,\infty ]$ with $\liminf_{t\to 0^+}h(t)/t=a$ and $\limsup_{t\to 0^+}h(t)/t=b$, where $0\le a\le b\le \infty$. Then, for every $x \in\ ]0, \infty [$ we have $h_*(x)\le ax$ and $h^*(x)\ge bx$.
\et

\pr Let $x>0$. By definitions of $h_*$ and $h^*$, for every positive integer $n$ we have $h_*(x)\le nh(x/n) \le h^*(x)$, that is,
\begin{equation}\label{eins}
h_*(x) \le x\cdot \frac{h(\frac{x}{n})}{\frac{x}{n}} \le h^*(x)\ .
\end{equation}
Since $h$ is increasing, for every $t$ such that $\frac{x}{n+1}\le t\le \frac{x}{n}$ we have
\[  \frac{h(\frac{x}{n+1})}{\frac{x}{n}}\le \frac{h(t)}{t} \le \frac{h(\frac{x}{n})}{\frac{x}{n+1}} \ .\]
Applying the limits inferior and superior to these inequalities as $t\to 0^+$ and $n\to \infty$ (with $\frac{n+1}{n}\to 1$) shows that
\begin{equation}\label{zwei}
\liminf_{n\to\infty} \frac{h(\frac{x}{n})}{\frac{x}{n}} \leq \liminf_{t\to 0^+} \frac{h(t)}{t}\ \ \ {\rm and} \ \ \
\limsup_{n\to\infty} \frac{h(\frac{x}{n})}{\frac{x}{n}} \geq \limsup_{t\to 0^+} \frac{h(t)}{t}\ .
\end{equation}
Combining (\ref{eins}) with (\ref{zwei}) now gives
\[ h_*(x) \le x\cdot \liminf_{t\to 0^+} \frac{h(t)}{t} = ax\ \ \ {\rm and} \ \ \
h^*(x) \ge x\cdot \limsup_{t\to 0^+} \frac{h(t)}{t} =bx\]
for every $x>0$, which completes the proof. \ebox

The values of $\liminf_{t\to 0^+} \frac{h(t)}{t}$ and $\limsup_{t\to 0^+} \frac{h(t)}{t}$ can be interpreted as the `lower' and `upper' slope of $h$ at the point $x=0$. The previous result may therefore be interpreted by saying that the values of $h_*$ and $h^*$ are to a large extent influenced by the values of the lower and upper slopes of $h$ at $0$.
\smallskip

\bc\label{con} Suppose that $h$ is an unary aggregation function on $[0, \infty]$ such that the derivative $h^{'}(0^+)$ exists and is equal to $c \in [0, \infty[$.
\smallskip

{\rm (1)}  If $h$ is convex on $[ 0, \infty[$, then $h_*(x)=cx$ and $h^*(x)=h(x)$ for every $x \geq 0$.
\smallskip

{\rm (2)}  If $h$ is concave on $[ 0, \infty[$, then $h_*(x)=h(x)$ and $h^*(x)=cx$ for every $x > 0$.
\ec

\pr For (1) it is sufficient to realize that $h(x)\geq ax$ for every $x \geq 0$, the claim then follows from Theorem \ref{ab} regarding $h_*$ and from \cite{Greco} regarding $h^*$. The proof of (2) is similar and therefore omitted. \ebox

\bc\label{ab1} For any real $a,b$ such that $0 < a < b < \infty$ there is an infinite number of smooth unary aggregation functions $h$ on $[0,\infty]$ such that $h_*(x)=ax$ and $h^*(x)=bx$ for each $x\in [0, \infty[$.
\ec

\begin{figure*}[h!]
 \centering
\includegraphics[width=0.6\textwidth]{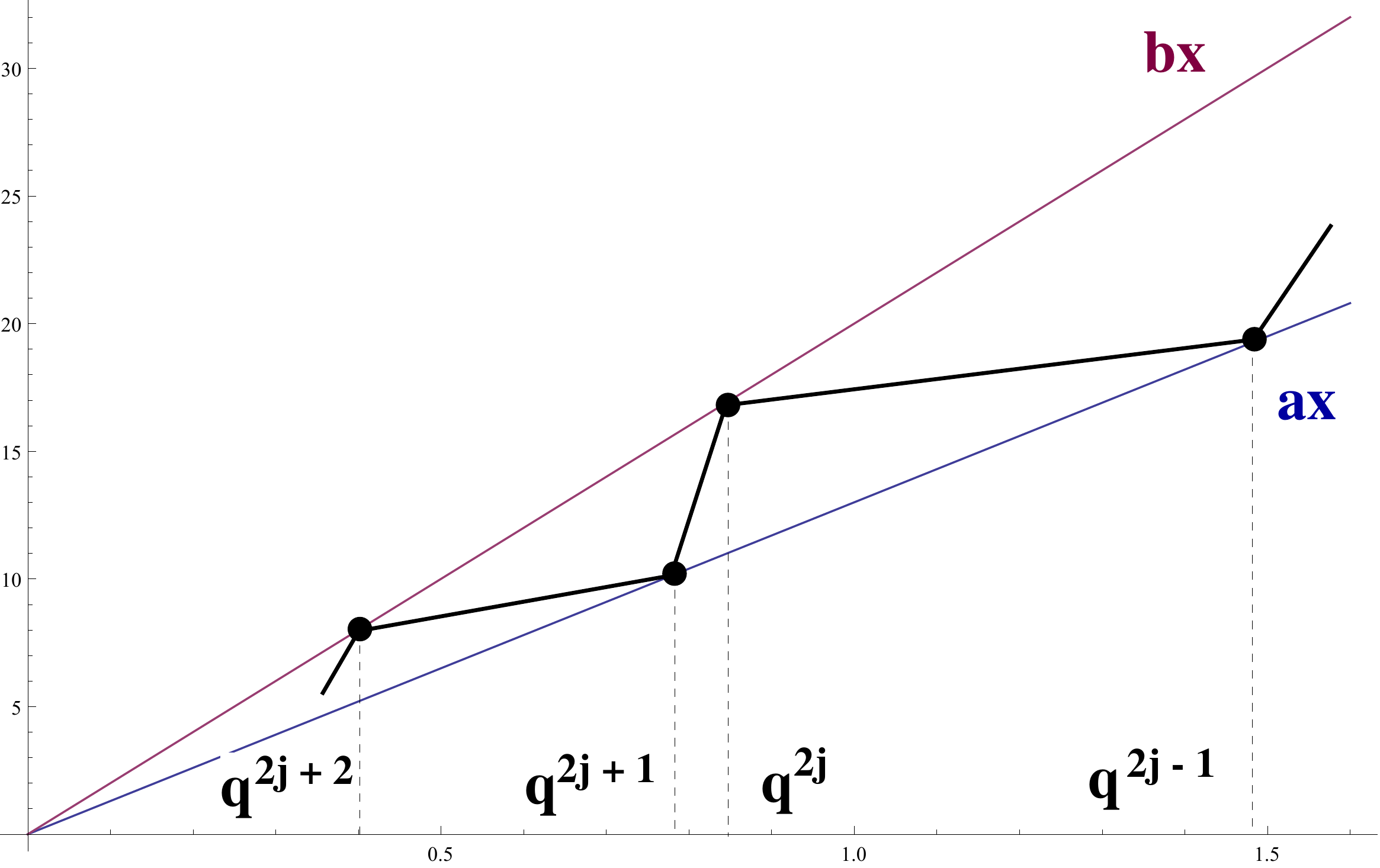}
\caption{A schematic drawing of a function h from the proof of Corollary 2.}
\end{figure*}

\pr Let $q$ be a positive real number such that $q < a/b <1$; note that $bq^{2j} < aq^{2j-1}$ for every positive integer $j$. Results from calculus now imply the existence of infinitely many smooth increasing functions $h(x)$ defined on $[ 0,\infty[$ such that $ax\le h(x)\le bx$ for every $x\in [ 0,+\infty[$, $h(q^{2j-1})=aq^{2j-1}$ and $h(q^{2j})=bq^{2j}$ for every positive integer $j$.
\smallskip

Since $ax\le h(x)\le bx$ for $x\in [ 0,\infty[$, we obviously have $ax\le h_*(x)$ and $h^*(x)\le bx$ for every $x\ge 0$. But we also have $\liminf_{t\to 0^+}h(t)/t=a$ and $\limsup_{t\to 0^+}h(t)/t=b$, because of the values of $h$ at points in the sequences $(q^{2j-1})_{j=1}^{\infty}$ and $(q^{2j})_{j=1}^{\infty}$, respectively.  By Theorem \ref{ab} we have $h_*(x)\le ax$ and $h^*(x)\ge bx$ for each $x\ge 0$, completing the proof. \ebox

Observe that the functions $h$ from Corollary \ref{ab1} have the property that $(h_*)^*(x)=ax < bx=(h^*)_*(x)$ for all $x>0$.
\smallskip

\bc\label{ab2} There is an infinite number of smooth aggregation functions $h$ on $[ 0,\infty]$ such that $h_*(x)=0$ for every $x< \infty$ and $h^*(x)=\infty$ for every $x>0$.
\ec

 \begin{figure*}[h!]
 \centering
 \includegraphics[width=0.6\textwidth]{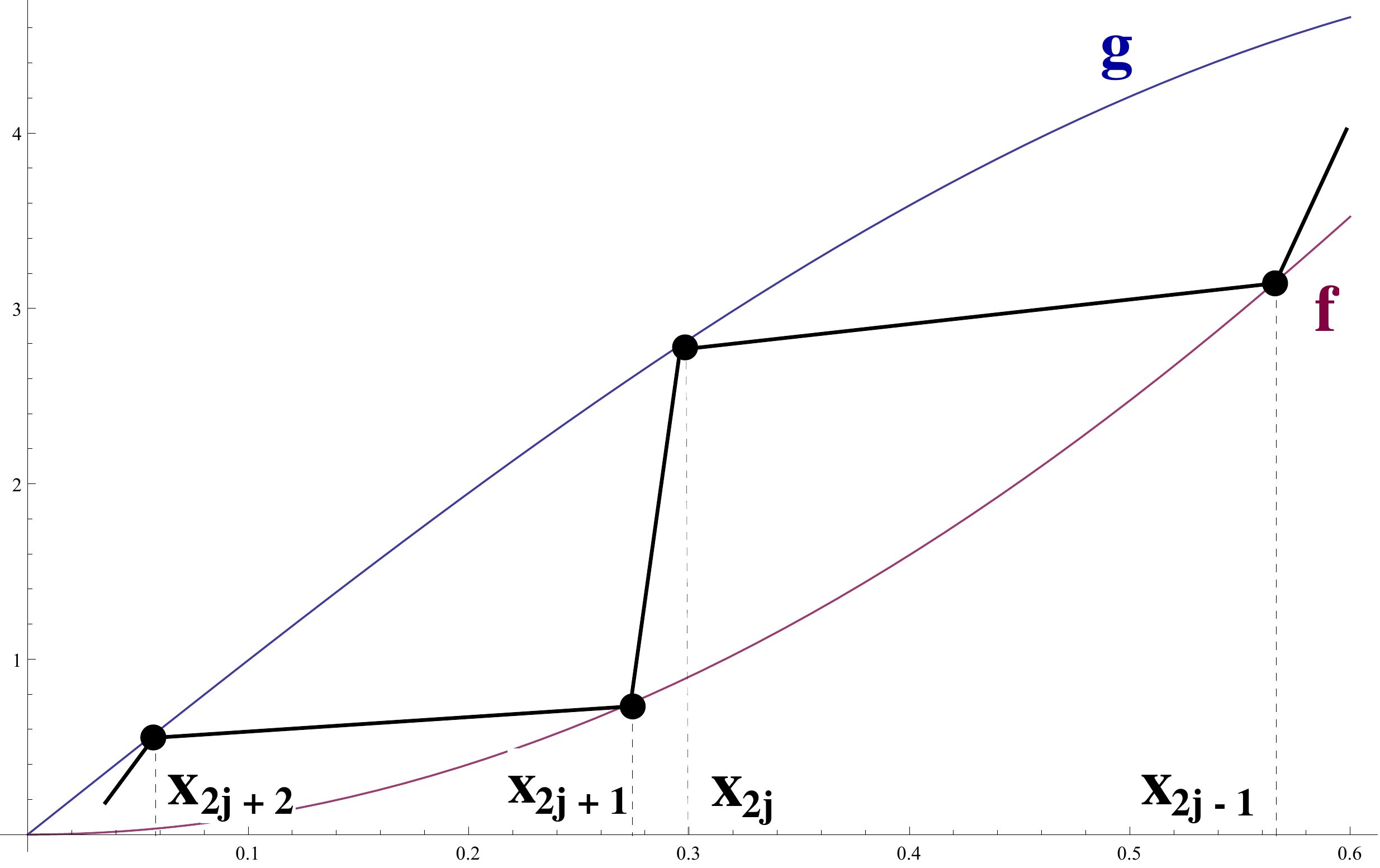}
  \caption{A schematic drawing of a function h from the proof of Corollary 3.}
  \end{figure*}

\pr For every positive integer $k$ let $x_k={2^{-2}}^k$. For $x\ge 0$ let $f(x)=x^{5/4}$ and $g(x)=x^{3/4}$; a straightforward calculation shows that $g(x_{2j}) < f(x_{2j-1})$ for every positive integer $j$. By known results from calculus there exists an infinite number of smooth increasing functions $h$ on $[ 0,\infty[$ such that $h(x_{2j-1})=f(x_{2j-1})$ and $h(x_{2j})=g(x_{2j})$ for all positive integers $j$.
\smallskip

Since for our function $h$ we have $\liminf_{t\to 0^+}h(t)/t=0$ and $\limsup_{t\to 0^+}h(t)/t=+\infty$ due to the values of $h$ at points in the sequences $(x_{2j-1})_{j=1}^{\infty}$ and $(x_{2j})_{j=1}^{\infty}$, the result is again a consequence of Theorem \ref{ab}. \ebox

The functions $h$ from Corollary \ref{ab2} have even a more striking property that $(h_*)^*(x)=0$ for every $x\ge 0$ while $(h^*)_*(x)=\infty$ for all $x>0$.
\smallskip

To conclude this section we underscore the fundamental role of Theorem \ref{ab} by the following complete characterization of one-dimensional degeneracies.

\bt\label{AB}  Let $h$ be a one-dimensional aggregation function on $[0,\infty[$. The following conditions are equivalent:
\smallskip

{\rm (a)} There exists an $x>0$ for which $h^*(x)=\infty$.
\smallskip

{\rm (b)} $h^*(x)=\infty$ for every $x>0$.
\smallskip

{\rm (c)} $\limsup_{t\to 0^+}\frac{h(t)}{t}=\infty$.
\smallskip

{\rm (d)} $\sup \{ \frac{h(t)}{t}\ |\ t\in]0,x]\}=\infty$ for some $x>0$.
\medskip

\noindent Similarly, the following statements are equivalent:
\medskip

{\rm (a')} There exists an $x>0$ for which $h_*(x)=0$.
\smallskip

{\rm (b')} $h_*(x)=0$ for every $x\ge 0$.
\smallskip

{\rm (c')} $\liminf_{t\to 0^+}\frac{h(t)}{t}=0$.
\smallskip

{\rm (d')} $\inf \{ \frac{h(t)}{t}\ |\ t\in]0,x]\}=0$ for some $x>0$.
   \et
\bigskip

\pr Clearly, the statements (c) and (d) are equivalent, and so are (c') and (d'). By Theorem 1, (c) implies (b) and (c') implies (b'). Trivially, (b) implies (a) and (b') implies (a'), and so one only has to prove that (a) implies (d) and (a') implies (d').
\smallskip

To show that (a) implies (d), we prove the contrapositive. Assume that $\sup \{ \frac{h(t)}{t}\ |\ t\in]0,x]\}=b < \infty$ for every $x>0$. This means that $h(t)\le bt$ for every $t\in]0,x]$. Thus, for every $n$-tuple $(x_1,x_2,\ldots,x_n)$ of positive real numbers such that $\sum_{i=1}^nx_i\le x$ we have $\sum_{i=1}^n h(x_i)\le b\sum_{i=1}^nx_i\le bx$. It follows that $h^*(x)\le bx < \infty$ for every $x>0$. \smallskip

Similarly, to show that (a') implies (d') we again proceed by proving the contrapositive. Suppose that $\inf \{ \frac{h(t)}{t}\ |\ t\in]0,x]\}=b >0$ for every $x>0$. This means that $h(t)\ge bt$ for every $t\in]0,x]$. Thus, for every $n$-tuple $(x_1,x_2,\ldots,x_n)$ of real numbers from $]0,x]$ such that $\sum_{i=1}^nx_i\ge x$ we have $\sum_{i=1}^n h(x_i)\ge b\sum_{i=1}^nx_i\ge bx$. It follows that $h^*(x)\ge bx > 0$ for every $x>0$. This completes the proof. \ebox

\section{The multidimensional case}\label{s3}

Based on the results in the one-dimensional case proved in Section \ref{s2} we continue by exhibiting examples of aggregation functions $A$ with the property that the values of $(A_*)^*$ are smaller than the values of $(A^*)_*$ for all non-zero vectors ${\bf x}\in [0,\infty[^n$. The method is based on the observation that the values of $A_*$ and $A^*$ on the entire space $[0,\infty[^n$ are influenced by the behavior of the one-dimensional diagonal function $A(x,\ldots,x)$ in an arbitrarily small neighbourhood of zero. We present details only regarding extensions of Corollary \ref{ab2} to arbitrary dimensions.
\smallskip

Let $\mu_{\bf x}$ be the arithmetic mean of the entries in ${\bf x}$, that is, if ${\bf x}=(x_1,x_2,\ldots, x_m)$, then $\mu_{\bf x}=(x_1+x_2+\ldots +x_m)/m$.
\smallskip

\bt\label{inf} There are infinitely many aggregation functions $A$ defined on $[0,\infty]^n$ such that $(A_*)^*({\bf x})=0$ for all ${\bf x}\in [0,\infty]^n$ while $(A^*)_*({\bf x})=\infty$ for all non-zero ${\bf x}\in [0,\infty[^n$.
\et

\pr It is sufficient to take any function $h$ from Corollary \ref{ab2} and define $A$ for every ${\bf x} \in [0,\infty[^n$ by $A({\bf x})=h(\mu_{\bf x})$. Letting $x=\mu_{\bf x}$, for any non-zero ${\bf x}$ we have, exactly as in the proof of Theorem \ref{ab},
\[ A_*({\bf x}) \le nh(x/n)\le x\cdot \liminf_{t\to 0^+} \frac{h(t)}{t}\ , \ \ {\rm and}\]
\[ A^*({\bf x}) \ge nh(x/n) \ge x\cdot \limsup_{t\to 0^+} \frac{h(t)}{t}\ .\]
The claim now follows from Corollary \ref{ab2}. \ebox

Clearly, one can use numerous other compositions of functions $h$ from Corollary \ref{ab2} with simple aggregation functions (such as weighted average, geometric means, etc.) to provide examples for Theorem \ref{inf}.
\smallskip

We also prove different sufficient conditions for the values of $A_*$ and $A^*$ to exhibit the extreme behavior described in Theorem \ref{inf}.
\smallskip

\bp\label{ser} Let $A$ be an aggregation function on $[0,\infty]^n$ and let $\tilde A(x)=A(x,x,\ldots,x)$ for every $x\ge 0$.
\smallskip

{\rm (1)} If there exists a divergent series $\sum_{j=1}^{\infty}a_j$ with decreasing positive terms such that the series $\sum_{j=1}^{\infty} \tilde A(a_j)$ converges, then $A_*({\bf x})=0$ for every ${\bf x}\in [0,\infty[^n$.
\smallskip

{\rm (2)} If there exists a convergent series $\sum_{j=1}^{\infty}a_j$ with decreasing positive terms such that the series $\sum_{j=1}^{\infty} \tilde A(a_j)$ diverges, then $A^*({\bf x})=+\infty$ for every non-zero ${\bf x}\in [0,\infty]^n$.
\ep

\pr Let $h(x)=A(x,x,\ldots,x)$ for $x\ge 0$. We show that the assumption of (1) implies that
$\lim\inf_{j\to \infty}h(a_j)/a_j=0$ and (2) implies that $\limsup_{j\to \infty}h(a_j)/a_j=+\infty$.
\smallskip

Indeed, suppose that $\lim\inf_{j\to \infty}h(a_j)/a_j=c>0$. This means that for every $\epsilon > 0$ we have $h(a_j) \geq (c - \epsilon)a_j$ for all but a finite number of positive integers $j$. But then, divergence of $\sum_j a_j$ would imply divergence of $\sum_j h(a_j)$, contrary to the assumption of (1). Similarly, if $\limsup_{j\to \infty}h(a_j)/a_j=c < +\infty$, then for every $\epsilon > 0$ we would have $h(a_j) \leq (c + \epsilon)a_j$ for all but finitely many $j$'s. But then convergence of $\sum_j a_j$ would imply convergence of $\sum_j h(a_j)$, a contradiction.
\smallskip

This shows that $\lim\inf_{t\to 0^+}h(t)/t=0$ in the case (1) and $\limsup_{t\to 0^+}h(t)/t=+\infty$ in the case (2). The result now follows from Theorems \ref{ab} and \ref{inf}. \ebox

As examples, one can take the function $f(x)=x^{1+\lambda}$ for some arbitrarily small $\lambda>0$, or $f(x)=x/\ln^2(x)$, and take for $\sum_{j=1}^{\infty}a_j$ the harmonic series to construct aggregation functions $A$ such that $\tilde A(x)=f(x)$ on an arbitrarily small interval $(0,\delta)$; the result of (1) then gives $A_*({\bf x})=0$ on $[0,\infty[^n$. Similarly, one can take the function $g(x)=x^{1-\lambda}$ for an arbitrarily small $\lambda>0$, or $g(x)=x\ln^2(x)$, and take for $\sum_{j=1}^{\infty}a_j$ the series with $a_j=g^{-1}(1/j)$ to obtain aggregation functions $A$ such that $\tilde A(x)=g(x)$ on an arbitrarily small interval $(0,\delta)$; the result of (2) then shows that $A^*({\bf x})=+\infty$ for all non-zero ${\bf x}\in [0,\infty[^n$.
\smallskip

With the help of the results of Section \ref{s2} we can also decide membership in ${\cal K}_n^*$ and ${\cal K}_*^n$ by looking at the one-dimensional case. For an aggregation function $A$ on $[0,\infty]^n$ let $\tilde A$ be defined on $[0,\infty]$ by letting $\tilde A(x)=A(x,x,\ldots,x)$ and, for every $i\in \{1,2,\ldots,n\}$ let $A_i$ be defined on $[0,\infty]$ by $A_i(x)=A(x{\bf e}_i)$, where ${\bf e}_i$ is the $i$-th unit vector.

\bt\label{K}
Let $A$ be an aggregation function on $[0,\infty]^n$. Then,
\bi
\item[{\rm (i)}] $A|[0,\infty[^n \in {\cal K}_n^*$ if and only if $\tilde A|[0,\infty[ \in {\cal K}_1^*$, and
\item[{\rm (ii)}] $A|[0,\infty[^n \in {\cal K}_*^n$ if and only if $A_i|[0,\infty[ \in {\cal K}_*^1$ for some $i\in \{1,2,\ldots,n\}$.
\ei.  \et

\pr
(i): We only need to show that $A^*({\bf x})=\infty$ for some non-zero vector ${\bf x}\in [0,\infty[^n$ if and only if $(\tilde A)^*(x)=\infty$ for some $x>0$. For the direct implication, by monotonicity it is sufficient to take, for a given non-zero vector ${\bf x}$, the value $x$ equal to the maximum of the coordinates of ${\bf x}$; the reverse implication follows by taking, for a given $x>0$, the vector ${\bf x}=(x,x,\ldots,x)$.
\smallskip

(ii): Here we only need to show that $A_*({\bf x})=0$ for some non-zero vector ${\bf x}\in [0,\infty[^n$ if and only if $(A_i)_*(x)=0$ for some $x>0$ and some $i\in \{1,2,\ldots,n\}$. For the direct implication it is sufficient to take an $i\in \{1,2,\ldots,n\}$ for which the $i$-th coordinate of ${\bf x}$ has a non-zero value $x$; by monotonicity we have $(A_i)_*(x)=0$. The reverse implication follows by simply taking ${\bf x}=x{\bf e}_i$.
\ebox

As an example of this result consider the aggregation function $A$ defined for any ${\bf x}\in [0,\infty]^3$ by taking $A({\bf x})$ to be the median of the coordinates of the vector ${\bf x}$. It follows immediately from Theorem \ref{K} that $A|[0,\infty[^3 \in {\cal K_3^*}$ but it does not belong to ${\cal K_*^3}$.
\smallskip

\section{Comparing $(A_*)^*$ with $(A^*)_*$}\label{s4}

In this section we prove that the values of $(A_*)^*$ can never be larger than the values of $(A^*)_*$ for an aggregation function $A$ defined on $[0,\infty]^n$. The proof strongly depends on the fact that $A$ is defined for non-zero vectors of $[0,\infty]^n$ that are arbitrarily close to the zero vector.
\smallskip

\bt \label{main} Let $A:\ [0,\infty]^n \to [0,\infty]$ be an aggregation function. Then $(A_*)^* \leq (A^*)_*$.  \et

\pr For our given aggregation function $A$ and for every $i \in \{1,\ldots,n\}$ let $A_i:\  [0, \infty[ \to [0, \infty[$ be the marginal function of $A$ considered before the proof of Theorem \ref{K}, given by $A_i(x)=A(x{\bf e}_i) = A(0, ... , x, ..., 0)$, where $x$ appears only in the $i$-th coordinate. It is not difficult to check that
\[ A_*(x_1, ... , x_n) \leq \sum \limits_{i=1}^n (A_i)_*(x_i) \ \ \ {\rm and} \ \ \ A^*(x_1, ... , x_n) \geq \sum \limits_{i=1}^n (A_i)^*(x_i)\ .\]
It follows that
\[ (A_*)^*(x_1, ... , x_n) \leq \sum \limits_{i=1}^n ((A_i)_*)^*(x_i) \ \ \ {\rm and} \ \ \ (A^*)_*(x_1, ... , x_n) \geq \sum \limits_{i=1}^n ((A_i)^*)_*(x_i)\ .\]
Let now
\[ a_i=\liminf_{t\to 0^+} \frac{A_i(t)}{t} \ \ \  {\rm and} \ \ \ b_i=\limsup_{t\to 0^+} \frac{A_i(t)}{t}  \ .\]
By Theorem \ref{ab} from Section \ref{s2} we have
\[ (A_i)_*(x) \leq a_ix \leq b_ix \leq (A_i)^*(x) \]
which implies that
\[((A_i)_*)^*(x) \leq a_ix \leq b_ix \leq ((A_i)^*)_*(x) \ .\]
Consequently, for every ${\bf x}=(x_1,x_2,\ldots,x_n)\in [0,\infty[^n$ we have
\[(A_*)^*({\bf x}) \leq \sum \limits_{i=1}^n ((A_i)_*)^*(x_i) \leq \sum \limits_{i=1}^n a_ix_i \leq \sum \limits_{i=1}^n b_ix_i \leq \sum \limits_{i=1}^n ((A_i)^*)_*(x_i)= (A^*)_*({\bf x}) \ ,\]
which completes the proof.  \ebox

We also include an observation about the other extreme, that is, when $(A_*)^*=(A^*)_*$. A characterization of aggregation functions $A$ on $[0,\infty]^n$ for which this equality holds is still an open problem. Examples of such non-linear functions can be constructed, for instance, using Corollary \ref{con} in conjunction with Theorem \ref{inf}. However, if $(A_*)^*=(A^*)_*$ we can at least say that the result is an additive function, i.e., a function $D$ satisfying $D({\bf x}+{\bf y})=D({\bf x}) + D({\bf y})$ for every ${\bf x},{\bf y}\in [0,\infty[^n$. We formulate and prove the result in a slightly larger generality.
\smallskip

\bc\label{last}
Let $B$ and $C$ be aggregation functions on $[0,\infty]^n$ such that $B_*({\bf x})=C^*({\bf x})$ for every ${\bf x} \in [0,\infty[^n$. Then, the function $D=B_*=C^*$ is additive, $D({\bf x})=\sum \limits_{i=1}^n w_ix_i$ for some $w=(w_1, ... , w_n)\in [0, \infty]^n$.
\ec

\pr By the results of \cite{Greco}, $B_*$ is sub-additive and $C^*$ is super-additive; it follows that $D$ is both sub- and super-additive and hence additive. \ebox

\bexam

{\rm (i):} Define $A\in {\cal A}_n$ by $A(x_1, \ldots , x_n)=max(x_1, \ldots , x_n)$. Then $A \in {\cal P}_n$ and $A=A_*$, $A^*(x_1, \ldots , x_n)=\sum \limits_{i=1}^n x_i$, and $(A_*)^*=(A^*)_*=A^*$.
\smallskip

{\rm (ii):} Consider $B, C \in {\cal A}_n$ given by $B(x_1, \ldots , x_n)=\ln(\Pi_{i=1}^n (1+x_i))$ and $C(x_1, \ldots , x_n)=(\sum\limits_{i=1}^n e^{x_i})-n$. Then $B_*=B$, $C^*=C$ and $B^*=C_*$ is given by $B^*(x_1, \ldots , x_n)=\sum \limits_{i=1}^n x_i$.

\eexam

\section{Concluding remarks}

Our main aim was to characterize proper aggregation functions $A$ with the property that $A_*$ (and, similarly, $A^*$) is proper as well. In the course of our investigation we obtained a number of related results on the behavior of the subadditive and superadditive transformations $ A_*$ and $A^*$. Regarding iterations of these transformations we proved that  $( A_*)^* \leq (A^*)_*$ and that these two functions may be arbitrarily far from each other.
\smallskip

In a nutshell and in a somewhat more abstract setting, we clarified constraints for aggregation functions admitting superadditive/subadditive transformations without a complete loss of information. The transformations considered are related to optimization problems, for example to production functions in economical problems. For a deeper discussion we recommend ${\cite{Greco}}$.
\bigskip

\noindent{\bf Acknowledgement.}~~ The author acknowledges support from the projects VEGA 1/0420/15 and APVV 0013/14.

\end{document}